\documentclass [letterpaper, 10pt, oneside] {article}

\usepackage{paralist}   
\usepackage{booktabs}   

\usepackage{graphicx}
\usepackage{tikz}
\usepackage{epstopdf}
\usepackage{amsfonts}
\usepackage{graphicx}
\usepackage{amsmath}
\usepackage{amsthm}
\usepackage{amssymb}
\usepackage{fullpage}
\usepackage{verbatim}
\usepackage{tikz,pgfplots}
\pgfplotsset{compat=1.8}
\usepackage{curves}
\usepackage{color}
\usepackage{url}
\usepackage[vlined,linesnumbered,algo2e]{algorithm2e}

\renewcommand{\O}{\mathcal{O}}
\renewcommand{\d}{\text{d}}

\renewcommand{\i}{\imath}
\theoremstyle{plain}
\newtheorem{thm}{Theorem}
\newtheorem{lemma}{Lemma}
\newtheorem{cor}{Corollary}
\newtheorem{defn}{Definition}

\usepgfplotslibrary{external}

\renewcommand{\d}{{\mathrm d}}

\begin{document}

\title{Gaussian Quadrature Rule using $\epsilon$-Quasiorthogonality}

\author{Pierre-David L\'etourneau,$^a$ Eric Darve,$^b$\\
$^a$ Institute of Computational and Mathematical Engineering (ICME), Stanford University, CA\\
$^b$ Mechanical Engineering Department, Stanford University, CA}

\date{}

\begin{abstract}
We introduce a new type of quadrature, known as approximate Gaussian quadrature (AGQ) rules using $\epsilon$-quasiorthogonality, for the approximation of integrals of the form $\int f(x) \, \d \alpha(x)$. The measure $\alpha(\cdot )$ can be arbitrary as long as it possesses finite moments $\mu_n$ for sufficiently large $n$. The weights and nodes associated with the quadrature can be computed in low complexity and their count is inferior to that required by classical quadratures at fixed accuracy on some families of integrands. Furthermore, we show how AGQ can be used to discretize the Fourier transform with few points in order to obtain short exponential representations of functions.
\end{abstract}


\maketitle

\section{Introduction}
\label{Intro}

In this paper, we present a new kind of quadrature rule for approximating integrals by sums of the form,
\begin{equation}
\label{discrete}
 \int f(x) \, \d \alpha(x) \approx \sum_{i=1}^n w_i f(x_i)
\end{equation}
having the following characteristics:
\begin{enumerate}
\item The measure $\alpha( \cdot )$ can be \emph{arbitrary} (positive, signed, complex, ...) as long as it satisfies some weak condition.
\item The nodes and weights associated with the quadrature rule can be obtained in low computational complexity through a simple numerical algorithm.
\item The quadrature is at least as accurate as the Gaussian quadrature rule, and in many cases is significantly more accurate.
\item Low-order rules are able to integrate high-order polynomials with high accuracy.
\end{enumerate}

The scheme presented in the current work uses a strategy similar to classical Gaussian quadrature rules (of which a few examples can be found in Table~\ref{classicalGaussQuad}). The Gaussian quadrature rule is designed to integrate exactly polynomials of degree at most $2n-1$ using $n$ quadrature points and weights:
\[ \int x^k \, \d \alpha(x) \]
for various weight functions $\frac{\d \alpha}{\d x}$ (see Table~\ref{classicalGaussQuad}).

\begin{table}[htbp]
\caption{Examples of classical Gaussian quadratures}
\begin{center}
\begin{tabular}{ccc} \toprule
Name & Interval & Measure ($\d \alpha / \d x$ ) \\ \midrule
Gauss-Legendre & $[-1,1]$  & $1$ \\
Gauss-Laguerre & $[0, \infty)$ & $e^{-x}$ \\
Gauss-Hermite & $(-\infty, \infty)$ & $e^{-x^2} $\\
Gauss-Jacobi & $(-1,1)$ & $(1-x)^{\alpha}(1+x)^{\beta} \; , \;\; \alpha, \beta > -1 $ \\ 
Chebyshev-Gauss (1st kind) & $(-1,1)$ & $1/\sqrt{1-x^2}$ \\
Chebyshev-Gauss (2nd kind) & $[-1,1]$ & $\sqrt{1-x^2}$ \\
\bottomrule
\end{tabular}
\end{center}
\label{classicalGaussQuad}
\end{table}

The paper is structured as follows. In Section~\ref{GQ}, a brief overview of classical Gaussian quadratures will be presented. In Section~\ref{sec:agq}, the concept of quasiorthogonal polynomial and approximate Gaussian quadrature will be introduced together with an error analysis. This will be followed in Section~\ref{NS} by numerical results. In the same section, we will discuss representations of functions by short sums of exponentials.

\section{Gaussian quadrature}
\label{GQ}

Gaussian quadratures are schemes used to approximate definite integrals of the form,
\begin{equation*}
\int_a^b f(x) \, \d \alpha(x)
\end{equation*}
by a finite weighted sum of the form,
\begin{equation*}
\sum_{n=0}^N w_n \, f(x_n)
\end{equation*}
where $a<b \in \mathbb{R}$. The coefficients $\{ w_n \}$ are generally referred to as the \emph{weights} of the quadrature, whereas the points $\{ x_n \}$ are referred to as the \emph{nodes}. An $(N+1)$-node Gaussian quadrature  can integrate polynomials up to degree $2N+1$ \emph{exactly} and is generally well-suited for the integration of functions that are well-approximated by polynomials.

In what follows, we will briefly describe how the nodes and weights of classical Gaussian quadratures can be obtained based on the classical theory of orthogonal polynomials. For this purpose, we shall denote the real and complex numbers by $\mathbb{R}$ and $\mathbb{C}$ respectively. $\alpha( \cdot )$ will represent an arbitrary measure (possibly complex) on $(\mathbb{R}, \mathcal{B})$ or $(\mathbb{C}, \mathcal{B})$ unless otherwise stated. Vectors are represented by lower case letter e.g., $v$. The $i^{th}$ component of a vector $v$ will be written as $v_i$, and we shall use super-indices of the form $v^{(j)}$ when multiple vectors are under consideration.

We begin by introducing four key objects: the orthogonal polynomials, the Lagrange interpolants, the moments of a measure $\alpha(\cdot)$ and the Hankel matrix associated with such a measure.

\begin{defn}{\bf (Orthogonal polynomial)}
A sequence $\{ p^{(k)} (x) \}_{k=0}^{\infty}$ of polynomials of degree $k$ is said to be a sequence of orthogonal polynomials with respect to a positive measure $\alpha(\cdot)$ if,
\begin{equation*}
\int p^{(k)}(x) p^{(l)} (x) \, \d \alpha(x)  =
\left\{
	\begin{array}{ll}
		0  & \mbox{if } k \not = l \\
		c_k & \mbox{if } k = l
	\end{array}
\right.
\end{equation*}
If in addition $c_k = 1 \; \forall k \in \mathbb{N}$, then the sequence is called \emph{orthonormal}.
\end{defn}
We shall hereafter assume that all such polynomials are monic, i.e., that they can be written as,
\begin{equation*}
p^{(k)} (x) = x^k + \sum_{n=0}^{k-1} p^{(k)}_n x^n
\end{equation*}
where $\{ p^{(k)}_n \}_{n=0}^{k-1}$ are some (potentially complex) coefficients. We then introduce Lagrange interpolants,

\begin{defn}{\bf (Lagrange interpolant)}
Given a set of $(d+1)$ data points $\{ (x_n, y_n) \}_{n=0}^d$, the Lagrange interpolant is the unique polynomial $L (x)$ of degree $d$ such that,
\begin{equation*}
L (x_n) = y_n , \; n = 0 ... d
\end{equation*}
It can be written explicitly as,
\begin{equation*}
L (x) = \sum_{n=0}^d y_n\, \ell_n (x)
\end{equation*}
where,
\begin{equation*}
\ell_n (x) = \prod_{\substack{{m=0}\\ {m \not = n}}}^d \frac{x-x_m}{x_n-x_m}
\end{equation*}
and $\ell_n (x)$ is referred to as the $n^{th}$ Lagrange basis polynomial.
\end{defn}
Finally we introduce the moments as well as the Hankel matrix associated with a measure $\alpha(\cdot)$,
\begin{defn}{\bf(Moment)}
Given an arbitrary measure $\alpha(\cdot)$ on $(\mathbb{R}, \mathcal{B})$, its $n^{th}$ moment $\mu_n$ is defined by the following Lebesgue integral,
\begin{equation*}
\mu_n = \int x^n \, \d \alpha(x)
\end{equation*}
whenever it exists.
\end{defn}

\begin{defn}{\bf (Hankel matrix)}
An $(N+1) \times (M+1)$ matrix $H$ is called the $(N+1) \times (M+1)$ Hankel matrix associated with the measure $\alpha( \cdot )$ if its entries take the form,
\begin{equation}
\begin{pmatrix}
\mu_0 & \mu_1 & \cdots & \mu_{M}  \\
\mu_1 & \mu_2 & \cdots & \mu_{M+1}\\
\vdots &\vdots &\vdots &\vdots \\
\mu_{N} & \mu_{N+1} & \cdots &  \mu_{N+M} \end{pmatrix}
\label{eq:hankel}
\end{equation}
i.e., $H_{ij} = \mu_{i+j}$, where $(\mu_0, \mu_1, \ldots, \mu_{N+M} )$ are the first $(N+M)$ moments of $\alpha(\cdot)$ whenever they exist.
\end{defn}

With these quantities we can now present the main results associated with classical Gaussian quadratures, 

\begin{thm}{\bf(Gaussian quadrature)}
Consider a positive measure $\alpha(\cdot)$ on $([a,b], \mathcal{B})$ (with $a,b \in \mathbb{R}$ potentially infinity) and a sequence of orthonormal polynomials $\{ p^{(k)} (x) \}_{k=0}^{\infty}$ with respect to $\alpha( \cdot)$. Then, the quadrature rule with nodes $\{ x_n\}_{n=0}^k$ consisting in the zeros of $p^{(k+1)}(x)$ and weights $\{ w_n \}_{n=0}^k$ given by,
\begin{equation*}
w_n = \int \ell_n (x) \, \d \alpha(x) 
\end{equation*}
integrates polynomials of degree $\leq 2k+1$ exactly.
\end{thm}

This is a classical result which can be found in \cite{Meurant} for instance. Explicit expression for the error incurred in the case of smooth integrand also exist.

To close this section, we introduce a further result characterizing the coefficients of the orthogonal polynomials $\{ p^{(k)}(x) \} $. As we shall see in the next section, this characterization lies at the heart of our scheme,
\begin{lemma}
Consider a positive measure $\alpha(\cdot)$ on $([a,b], \mathcal{B})$ (with $a,b \in \mathbb{R}$ potentially infinity) and a sequence of orthogonal polynomials $\{ p^{(k)} (x) \}_{k=0}^{\infty}$ with respect to $\alpha( \cdot)$. Then, the coefficients $\{ p^{(k+1)}_n \}_{n=0}^k$ of the $(k+1)^{th}$ orthogonal polynomial $p^{(k+1)}(x)$ satisfy the following Hankel system,
\[ Hp =  
\begin{pmatrix}
\mu_0 & \mu_1 & \cdots & \mu_{k+1}  \\
\mu_1 & \mu_2 & \cdots & \mu_{k+2}\\
\vdots &\vdots &\vdots &\vdots \\
\mu_{k} & \mu_{k+1} & \cdots &  \mu_{2k+1} 
\end{pmatrix}
\begin{pmatrix}
p_0^{(k+1)}  \\
p_1^{(k+1)}  \\
\cdots  \\
p_{k}^{(k+1)}
\end{pmatrix} = 0\]  
where $\{ \mu_n \}$ are the moments of the measure $\alpha ( \cdot ) $, whenever they exist.
\label{AGQ:characterization}
\end{lemma}

\begin{proof}
First write,
\begin{equation*}
p^{(k+1)} (x) = \sum_{n=0}^{k+1} p^{(k+1)}_n x^n
\end{equation*}
Let $0 \leq j \leq k $. Then, from orthogonality we have,
\begin{equation*}
0 = \int p^{(k+1)} (x) \, x^j \, \d \alpha(x) =  \sum_{n=0}^{k+1}  p^{(k+1)}_n  \int  x^{n+j}\, \d \alpha(x) =  \sum_{n=0}^{k+1}  p^{(k+1)}_n  \mu_{n+j}
\end{equation*}
Putting all these equations in matrix form provides the desired result.
\end{proof}

The Hankel matrices associated with positive measures commonly encountered with classical Gaussian quadratures have been the subject of extensive study in the past (known as the moment problem). In some cases, they can be proved to be invertible although extremely ill-conditioned (see $\cite{Shohat}$ for details). On the other hand, less is known per regards to more general measures. In any case, in the event where the resulting Hankel matrix would be invertible, it can be expected to be ill-conditioned. Indeed, as an example it can be shown that for a large class of positive measures, the smallest eigenvalue of the $N \times N$ associated Hankel matrix scales like $\O \left (\frac{ \sqrt{N} }{\sigma^{2N}} \right ) $, where $\sigma$ depends only on the interval considered and is equal to $(1+\sqrt{2})$ for the interval $[-1,1]$ (see \cite{Widom:1966}).

The question we treat in the next section is whether such Hankel matrices arising from arbitrary measures can be used to derive Gaussian-like quadratures, and what this inherent ill-conditioning entails.

\section{Approximate Gaussian quadrature (AGQ)}
\label{sec:agq}

In this section, we describe the concept of approximate Gaussian quadrature. For this purpose, we will need the concept of $\epsilon$-quasiorthogonal polynomial, which we introduce for the first time below. Before doing so however, we first point to the following key observation.

\bigskip

\begin{thm}
Let $H$ be a $N \times M$ with rank $0 < d < M $. Then, there exists $D \le d+1$ and a vector $a \not = 0$ such that
\begin{equation*}
Ha = 0, \quad \text{with $a_i = 0$ for all $i > D$}
\end{equation*}
\label{AGQ:low_rank}
\end{thm}

\begin{proof}
The rank of $H$ is $d$. Therefore if we consider the first $d+1$ columns for $H$ they are linearly dependent. Denote $D$ the smallest integer such that the first $D$ columns of $H$ are linearly dependent. We have $D \le d+1$ and, by definition, there is $a \neq 0$ such that $Ha = 0$ with $a_i = 0$, $i>D$.
\end{proof}

\bigskip

We also have the following corollary,
\begin{cor}
\label{AGQ:quasiortho_cor}
Assume that  the $N \times (N+1)$ Hankel matrix $H$ associated with the measure $\alpha(\cdot)$ exists. If $H$ has rank $d<N$ then there exists a nontrivial polynomial $p(x)$ with degree $(D-1)$ where $D \leq d+1$ such that,
\begin{equation*}
\int p(x) x^j \, \d \alpha(x) = 0 
\end{equation*}
for all $j = 0,$ \ldots, $N$.
\end{cor}

\begin{proof}
Let $K$ be such that
 \begin{equation*}
D = \inf \{  0\leq n \leq N : \mathrm{rank}( H(:, 1:n) )= n \} 
\end{equation*}
where $H(:, 1:n)$ is the matrix containing the first $n$ columns of $H$. By theorem \ref{AGQ:low_rank}, there exists a vector $a \not = 0$ such that $Ha = 0$ and $a_i = 0$ for $i>D$.\\

Let $p(x)$ be the polynomial with coefficients given by $a$, i.e.
\begin{equation*}
p(x) = \sum_{n=0}^D a_n x^n
\end{equation*}
Then,
\begin{align*}
\int p(x) x^j \, \d \alpha(x) &= \int  \sum_{n=0}^{D} a_n x^{n+j}  \, \d \alpha(x)\\
&= \sum_{n=0}^D a_n \mu_{n+j} \\
& = (H a)_j = 0
\end{align*}
since $a$ belongs to the null-space of $H$.
\end{proof}

\begin{figure}[htbp]
\begin{center}
 \begin{tikzpicture}[scale=0.9]
	\begin{axis}[
	ylabel={$\log_{10}$ of absolute error},
 	xlabel={Index of singular value},
		ymode=log,
		width=0.6\textwidth,
		height=200pt,
		xmin=0,
		xmax=50,
		ymin=1e-18,
		ymax=1e2,
		xtick={0, 5, 10, 15, 20, 25, 30, 35, 40, 45, 50},
		ytick={1e-18, 1e-16,1e-14, 1e-12, 1e-10,1e-08, 1e-06, 1e-4, 1e-2, 1e0, 1e2}]
	
	\input{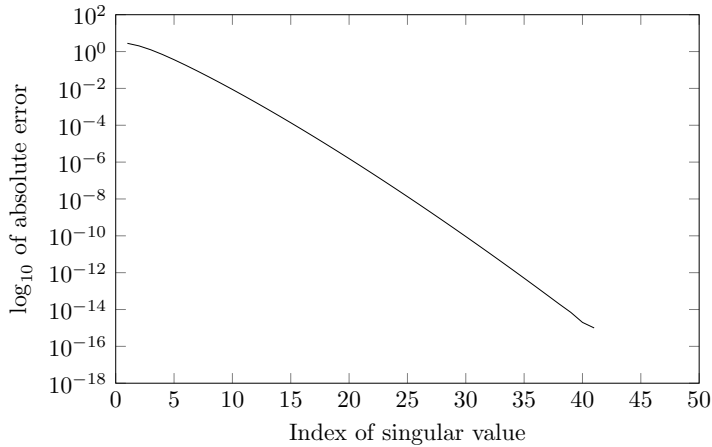}

	\end{axis}
\end{tikzpicture}%
\end{center}
\caption{$\log_{10}$ of singular values of Hankel matrix ($N = 250$) associated with the Lebesgue measure on $[-1,1]$.}
\label{svd}
\end{figure}

The consequences of this corollary are far-reaching and constitute the crux of the scheme presented here.  Indeed, although we do not generally expect the Hankel matrix $H$ associated with some measure $\alpha( \cdot)$ to be \emph{exactly} low-rank as in the case of Theorem \ref{AGQ:low_rank} (e.g., $H$ has full rank in the case of classical Gaussian quadratures) we can expect that in some cases $H$ will be \emph{approximately} low rank. In other words, given $0 < \epsilon  \ll 1 $ we expect,
\begin{equation*}
D \approx \max \{ 1 \leq i \leq N :  \sigma_i  > \epsilon \, \sigma_1 \}
\end{equation*}
where $\{ \sigma_i \}$ are the singular values of $H$, to be much smaller than $N$, i.e., $D \ll N$. We show for instance in Figure \ref{svd} the first $50$ singular values of the Hankel matrix ($N=250$) associated with the Lebesgue measure in $[-1,1]$. The $y$-axis scales as a logarithm in base $10$, and it is seen that the singular values decay faster than exponentially.

In light of the above discussion, we might expect in these circumstances the existence of a polynomial $p(x)$ of degree $D \approx \max \{ 1 \leq i \leq N :  \sigma_i  > \epsilon \sigma_1 \}$ such that 
\begin{equation*}
\left | \int p(x) x^j \, \d \alpha(x)  \right | \lesssim \epsilon
\end{equation*}
for all $0 \leq j \leq N$, and this leads us to the introduction of the concept of $\epsilon$-quasiorthogonal polynomial which we now define,
\begin{defn}
A polynomial $p(x)$ is called $\epsilon$-quasiorthogonal of order $N$ with respect to the measure $\alpha(\cdot)$ and the basis $\{ L_n (x) \}$ if,
\begin{equation*}
\left | \int p(x) \, L_n(x)  \, \d \alpha(x) \right |  \leq \epsilon
\end{equation*}
for all $n=0,$ \ldots, $N$.
\end{defn}
Importantly, this definition imposes no restriction per regards to the measure $\alpha(\cdot)$, in opposition with orthogonal polynomials which demand the measure to be positive (\cite{Szego}). In this sense, the relation described is not one of orthogonality for it is not possible to define a nondegenerate  inner-product unless $\alpha( \cdot )$ is positive. This is why we chose the name \emph{quasi}-orthogonal. We also note that given $\epsilon \geq \sigma_N(H)$ such polynomial always exists for it suffices to pick $a$ aligned with the right singular vector associated with the smallest singular value $\sigma_N(H)$.

From a computational standpoint, there exists an efficient scheme to find such polynomials given a measure $\alpha( \cdot)$ and some $\epsilon>0$. This is the subject of Section \ref{AGQ:comp}. For the remaining of this section, we will focus on demonstrating how such polynomials can be used to obtain efficient quadratures. As will be shown, the construction of the scheme shares a lot with that of classical Gaussian quadrature. This is what constitutes the origin of the denomination.

We will need the following technical lemma which proof is provided in appendix,
\begin{lemma}
\label{AGQ:tech_lemma}
Let $\alpha(\cdot)$ be an arbitrary measure on $(\mathbb{R}, \mathcal{B})$, and
\begin{equation*}
p(x) =x^{d+1} +  \sum_{n=0}^{d} p_n x^n
\end{equation*} 
 be a monic $\epsilon$-quasiorthogonal polynomial of degree $d+1$ and order $N$ associated with $\alpha(\cdot)$. Further, let,
\begin{equation*}
f(x) = \sum_{n=0}^{N+d} f_n  \, L_n(x) 
\end{equation*}
be some polynomial of degree $N+d$ and $ \tilde{q}(x) = \sum_{n=0}^{d} \tilde{q}_n x^n$ be the Lagrange interpolant of $q(x)$ associated with the zeros of $p(x)$. Finally, let $r(x) = \sum_{n=0}^{N} r_n x^n$ be the unique polynomial such that $ q(x) - \tilde{q}(x)= p(x) r(x) $. Then,
\begin{align*}
 \sum_{n=0}^N |r_n|  \leq \lVert\Gamma^{-1} \bar{q} \rVert_{1} 
\end{align*}
where $\bar{q} =\left( q_{d+1}, q_{d+2}, \ldots, q_{N+d}\right )^T $ and $\Gamma$ is the $N \times N$ Toeplitz matrix such that $[ \Gamma ]_{i,j} = p_{j-i}$ if $0\leq j-i \leq d$ and $0$ otherwise.
\label{AGQ:errorlemma}
\end{lemma}

We are now ready to prove our main theorem.

\begin{thm}{\bf(Approximate Gaussian quadrature)}
\label{AGQ:AGQ_thm}
Consider an arbitrary measure $\alpha(\cdot)$ on $(\mathbb{R}, \mathcal{B})$. Let $p(x)$ be a monic $\epsilon$-quasiorthogonal polynomial of degree $d+1$ and order $N$ with respect to $\alpha( \cdot )$, where $0<\epsilon<1$. Then, the quadrature rule with nodes $\{ x_n \}_{n=0}^{d}$ consisting in the zeros of $p(x)$ and weights $\{ w_n \}_{n=0}^{d}$ given by,
\begin{equation}
\label{AGQ:AGQ_thm:weight}
w_n = \int \ell_n (x) \, \d \alpha(x) 
\end{equation}
where $\ell_n (x)$ is the $n^{th}$Lagrange basis polynomial associated with the nodes, integrates polynomials $q(x)$ of degree $\leq N+d$ with an error bounded by,
\begin{equation*}
 \left |  \int  q(x) \, \d \alpha(x) -   \sum_{n=0}^d w_n \, q(x_n) \right |  \leq  \lVert\Gamma^{-1} \bar{q} \rVert_{1}  \, \epsilon
\end{equation*}
where $\{ q_n \}_{n=0}^{N+d}$ are the coefficients of $q(x)$, $\bar{q} =\left( q_{d+1} , q_{d+2}  , \ldots, q_{N+d}\right )^T $ and $\Gamma$ is the $N \times N$ Toeplitz matrix such that $[ \Gamma ]_{i,j} = p_{j-i}$ if $0\leq j-i \leq d$ and $0$ otherwise.
\end{thm}

\begin{proof}
Let $q(x)$ be a polynomial of degree $ N + d$ and consider the Lagrange interpolant at the nodes $\{ x_n \}$,
\begin{equation*}
\tilde{q} (x) = \sum_{n=0}^d  q(x_n) \ell_n (x)
\end{equation*}
Then consider,
\begin{align*}
I = \int \left [ q(x) - \tilde{q}(x) \right ] \, \d \alpha(x)
\end{align*}
The quantity $[q(x) - \tilde{q}(x)]$ is a polynomial of degree at most $(N+d)$ and has zeros located at each of the nodes $\{ x_n \}_{n=0}^d$. Therefore, by the factorization theorem for polynomials we can write,
\begin{equation*}
q(x) - \tilde{q}(x) = \prod_{n=0}^d (x - x_n) \, r(x)
\end{equation*}
where $r(x)$ is a polynomial of degree at most $N$. We further note that $\prod_{n=0}^d (x - x_n)$ is a monic polynomial of degree $d+1$ with zeros at $\{ x_n \}_{n=0}^d$ just as $p(x)$. Since monic polynomials are uniquely characterized by their roots we have,
\begin{equation*}
\prod_{n=0}^d (x - x_n) = p (x)
\end{equation*}
Therefore,
\begin{align*}
|I| &= \left | \int p(x) r(x) \, \d \alpha(x) \right | \leq \sum_{n=0}^N | r_n | \,\left |   \int  p(x) \, x^n \, \d \alpha(x) \right |  \leq  \sum_{n=0}^N | r_n | \, \epsilon \\
\end{align*}
where we used the $\epsilon$-quasiorthogonality of $p(x)$. Finally, thanks to Lemma \ref{AGQ:tech_lemma} we get,
\begin{equation*}
 \left |  \int  q(x) \, \d \alpha(x) -   \sum_{n=0}^d w_n \, q(x_n) \right |  \leq  \lVert \Gamma^{-1} \bar{q} \rVert_{1} \epsilon
\end{equation*}
\end{proof}

Interestingly, the above analysis reveals that an AGQ of order $d$ is in fact \emph{exact} for polynomials of degree $\leq d$.

Some advantages of AGQ is that there is no need for the measure $\alpha( \cdot )$ to have any specific properties beyond the existence of moments of high-enough order. Furthermore, the problem of the existence and uniqueness of the solution to the Hankel system is of no importance; in fact, the larger the null-space of $H$ the better it is. 

Both characteristics are in sharp contrast with common wisdom regarding classical Gaussian quadratures. First, the positivity of the measure is key in proving the existence of a sequence of orthogonal polynomials necessary to build a classical quadrature (see \cite{Meurant}, Theorem 2.7). Secondly, the notion of orthogonality is at the heart of modern numerical schemes used to obtain nodes and weights for it gives rise to a three-term recurrence relation that is thoroughly exploited computationally (see \cite{Golub:1969, Meurant}).

\subsection{Computational considerations}
\label{AGQ:comp}

The first computational issue we describe here is that of finding an adequate $\epsilon$-quasiorthogonal polynomials of order $N$ given a measure $\alpha(\cdot)$ on $(\mathbb{R}, \mathcal{B})$, some $N \in \mathbb{N}$ and some value $ 0< \epsilon$. For this purpose, we note that a sufficient condition for a monic polynomial $p(x)$ of degree $(d+1)$ to fall within this category is to satisfy the following inequality,
\begin{equation*}
\lVert H(N,d) \, \bar{p} + h(d) \rVert_{\infty} \leq \epsilon
\end{equation*}
where $\bar{p} = [p_0, \, p_1 \, , \ldots, p_{d}]^T$, $p(x) = x^{d+1} + \sum_{n=0}^d p_n x^n$, $ h(d) = [\mu_{d+1}, \, \mu_{d+2} \, , \ldots, \mu_{d+N+1}]^T$ and $H(N,d)$ is the $(N+1)\times (d+1)$ Hankel matrix associated with the measure, i.e.
\[ H(N,d)  =
\begin{pmatrix}
\mu_0 & \mu_1 & \cdots & \mu_{d}  \\
\mu_1 & \mu_2 & \cdots & \mu_{d+1}\\
\vdots &\vdots &\vdots &\vdots \\
\mu_{N} & \mu_{N+1} & \cdots &  \mu_{d+N}
\end{pmatrix}
\] 
The proof is analogous to that of Corollary \ref{AGQ:quasiortho_cor} and uses the definition of quasiorthogonal polynomials.

This inequality provides a constructive way for finding an $\epsilon$-quasiorthogonal polynomial of small degree. This is described in Algorithm \ref{poly_alg}; note that we replace the $\lVert \cdot \rVert_{\infty}$ norm by the more computationally-friendly $\lVert \cdot \rVert_2$ norm which is equivalent.

\begin{algorithm2e}[htbp]
Let $(d+1) =  \mathrm{rank}(H,\delta)$\\
Solve $ \min_p \lVert H(N,d) \, p + h(d) \rVert_2  $ \\
\While{$\lVert H(N,d) \, p + h(d) \rVert_2 > \epsilon $} {
	$d = d+1$ \\
	Solve $ \min_p \lVert H(N,d) \, p + h(d) \rVert_2  $ 
}
\caption{Pseudo-code to determine $\epsilon$-quasiorthogonal polynomial of order $N$ of small degree for desired accuracy $\delta$.}
\label{poly_alg}
\end{algorithm2e}

{\bf Note: } The quadrature obtained from $p(x)$ integrates polynomials of degree $N+d$ with error prescribed by Theorem \ref{AGQ:AGQ_thm}. This error term involves the norm of the inverse of a matrix $\Gamma$ which is upper-triangular, Toeplitz with diagonal entries all equal to $1$ and remaining entries depending on the coefficients of the polynomial $p(x)$. In order to guarantee that an AGQ integrates polynomials of degree $\leq N+d$ with accuracy $\delta$ say, it is sufficient to set $\epsilon \leq \frac{\delta}{C}$ and constrain $p(x)$ to be such that $\lVert \Gamma^{-1} \rVert_{\infty} \leq C$ for some $C>0$. Upon obtaining some characterization of the set $\mathcal{S}_C := \{ p(x) : \lVert \Gamma^{-1} \rVert_{\infty} \leq C \}$, one could potentially carry out the steps described in Algorithm \ref{poly_alg} while restraining the solution to $\mathcal{S}_C$. One would thus guarantee the accuracy of the AGQ \emph{a priori}. Unfortunately, such characterization is not readily available so one is left with the \emph{a posteriori} estimates of Theorem \ref{AGQ:AGQ_thm}. On the other hand, numerical experiments point to the fact that the product $\lVert \Gamma^{-1} \rVert_{\infty} \, \epsilon$ does indeed decay in a fast manner as a function of the degree of $p(x)$, for $\bar{p}$ the solution of the least-squares problem having the smallest norm in Algorithm \ref{poly_alg}. In short, although AGQ in its current state performs well, some improvements are still possible. This constitutes a topic for future research.

Once such polynomial has been obtained, its roots constitute the nodes of the approximate Gaussian quadrature as per Theorem \ref{AGQ:AGQ_thm}. The cost of solving a thin $(N+1) \times (d+1)$ least-squares problem is $\O( [N+1) + (d+1)/3] (d+1)^2 ) $ (see \cite{Golub}). Since in general we expect $d \ll N$ the cost is \emph{linear} in $N$. Also, each step of the while loop constitutes a rank-1 update of the system, so $p$ can be recomputed cheaply.

Another great computational aspect of the scheme is the availability of a simple analytical formula for the computation of the weights. Indeed, from Theorem \ref{AGQ:AGQ_thm} we have,
\begin{equation*}
w_n = \int \ell_n (x) \, \d \alpha(x) = \int \sum_{k=0}^d [\ell_n]_k x^k \, \d \alpha(x) = \sum_{k=0}^d [\ell_n]_k \, \mu_k
\end{equation*}
where $[\ell_n]_k $ is the $k^{th}$ coefficient of the $n^{th}$ Lagrange basis polynomial $ \ell_n (x)$, which can be obtained cheaply from the zeros of $ \ell_n (x)$, i.e., the nodes of the quadrature. We also noticed that it is generally possible to neglect nodes associated with small weights when such are present. This further reduces the cost of the method.

As a final comment, the accuracy of the scheme is highly dependent on the accuracy of the nodes. For this reason, we recommend performing the computations in extended arithmetic. In this paper, we used $Maple^\copyright$ in order to compute the nodes and weights of each approximate quadrature with high precision.

\section{Numerical simulations}
\label{NS}

In this section, we demonstrate the efficiency and the versatility of the scheme through a few numerical examples. In section \ref{NS:classical}, we compare fixed-order approximate Gaussian quadratures (AGQ) with two types of classical Gaussian quadratures (Gauss-Legendre and Gauss-Chebyshev) on monomials $x^n$ of increasing degree and show how it quickly becomes advantageous to use an approximate quadrature in those cases. Then in Section \ref{NS:sing}, we give examples related to functions with an integrable singularity at the origin.

In section \ref{NS:trig}, we show how the scheme can be applied to monomials on the complex circle, i.e., functions of the form $e^{\i n x}$ where $0 \leq n$. The resulting quadratures are then used in Section \ref{NS:Beylkin} to obtain approximations of functions through short exponential sums which is related to the method of Beylkin \& Monz\'on \cite{Beylkin:2005, Beylkin:2010}.

\subsection{Comparison with classical quadratures}
\label{NS:classical}
In this section, we compare results between the approximate Gaussian quadrature scheme, the Gauss-Legendre  $\left (  \d \alpha(x) = \d x  \right ) $ and Gauss-Chebyshev $\left (  \d \alpha(x) = \frac{1}{\sqrt{1-x^2} }\d x \right ) $ quadrature. 

\subsubsection{Integration of monomials}

For this benchmark, we fix the order ($N$ in Section \ref{AGQ:comp}) and study the error in approximating integrals of the form,
 \begin{equation*}
\int_{-1}^1 x^n \, \d \alpha (x)
\end{equation*}
through quadratures involving different number of nodes ($d$ in Section \ref{AGQ:comp}) where $n$ varies between $0$ to $700$.

\begin{figure}[htbp]
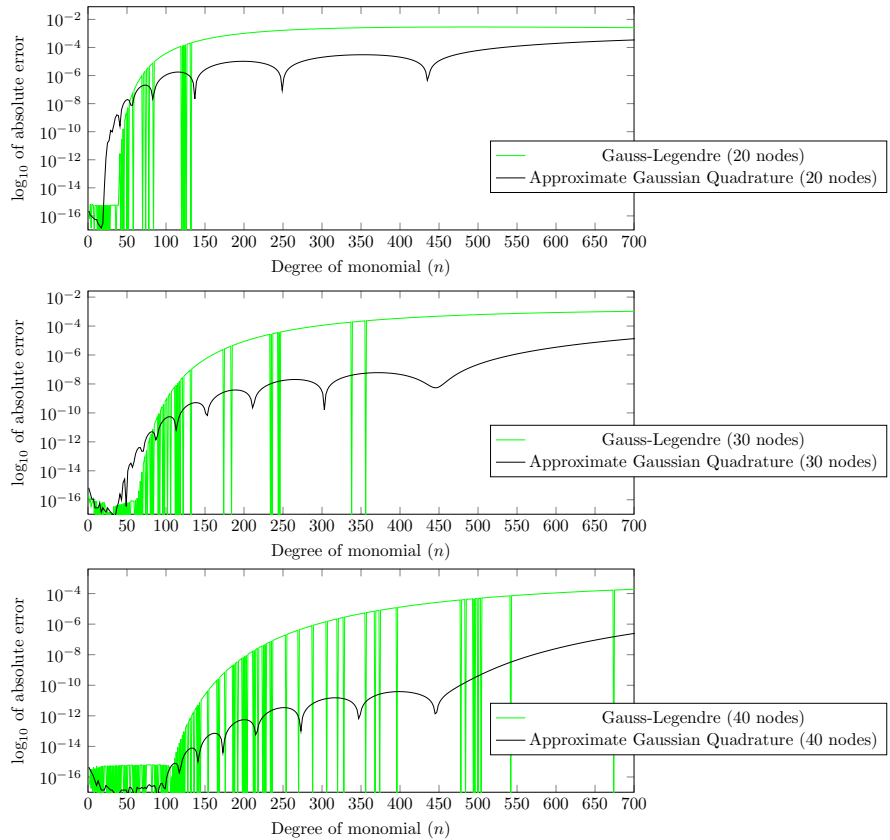

\begin{center}
 \begin{tikzpicture}[scale=0.65]
	\pgfplotsset{every axis legend/.append style={at={(1.1,0.4)},anchor=north}}	
	\begin{axis}[
	ylabel={$\log_{10}$ of absolute error},
 	xlabel={Degree of monomial ($n$)},
		ymode=log,
		width=0.75\textwidth,
		height=175pt,
		xmin=0,
		xmax=700,
		ymin=1e-17,
		ymax=0,
		xtick={0, 50,100,150,200,250,300,350, 400,450,500,550,600,650,700},
		ytick={1e-16,1e-14, 1e-12, 1e-10,1e-08, 1e-06, 1e-4, 1e-2}]
	
	\input{GaussLegendre20}

	\legend{Gauss-Legendre (20 nodes) \\%
	    Approximate Gaussian Quadrature (20 nodes) \\%
	}	
	\end{axis}
\end{tikzpicture}%

\begin{tikzpicture}[scale=0.65]
	\pgfplotsset{every axis legend/.append style={at={(1.1,0.4)},anchor=north}}	
 	\begin{axis}[
		ylabel={$\log_{10}$ of absolute error},
	 	xlabel={Degree of monomial ($n$)},
 		ymode=log,
 		width=0.75\textwidth,
 		height=175pt,
 		xmin=0,
 		xmax=700,
 		ymin=1e-17,
 		ymax=0,
 		xtick={0, 50,100,150,200,250,300,350, 400,450,500,550,600,650,700},
 		ytick={1e-16,1e-14, 1e-12, 1e-10,1e-08, 1e-06, 1e-4, 1e-2}]

		\input{GaussLegendre30}
 
	 	\legend{
	 	    Gauss-Legendre (30 nodes) \\%
		    Approximate Gaussian Quadrature (30 nodes) \\%
 		}
 	
 	\end{axis}
	\end{tikzpicture}%
	
\begin{tikzpicture}[scale=0.65]
	\pgfplotsset{every axis legend/.append style={at={(1.1,0.4)},anchor=north}}	
 	\begin{axis}[
		ylabel={$\log_{10}$ of absolute error},
 		xlabel={Degree of monomial ($n$)},
 		ymode=log,
 		width=0.75\textwidth,
 		height=175pt,
 		xmin=0,
 		xmax=700,
 		ymin=1e-17,
 		ymax=0,
 		xtick={0,50,100,150,200,250,300,350, 400,450,500,550,600, 650,700},
 		ytick={1e-16,1e-14, 1e-12, 1e-10,1e-08, 1e-06, 1e-4, 1e-2}]

		\input{GaussLegendre40}
 
 		\legend{
			     Gauss-Legendre (40 nodes) \\%
	     Approximate Gaussian Quadrature (40 nodes) \\%
 		}
 	
 	\end{axis}
	\end{tikzpicture}%
\end{center}
\caption{Comparison between the absolute error incurred in the evaluation of the integral $\int_{-1}^1 x^n \, \d x$ through an Approximate Gaussian Quadrature of order $N=350$ (black)  and a Gauss-Legendre quadrature (green) for different number of nodes. Top: 20 nodes, Middle: 30 nodes, Bottom: 40 nodes}
\label{GLcompare}
\end{figure}

\begin{figure}[htbp]
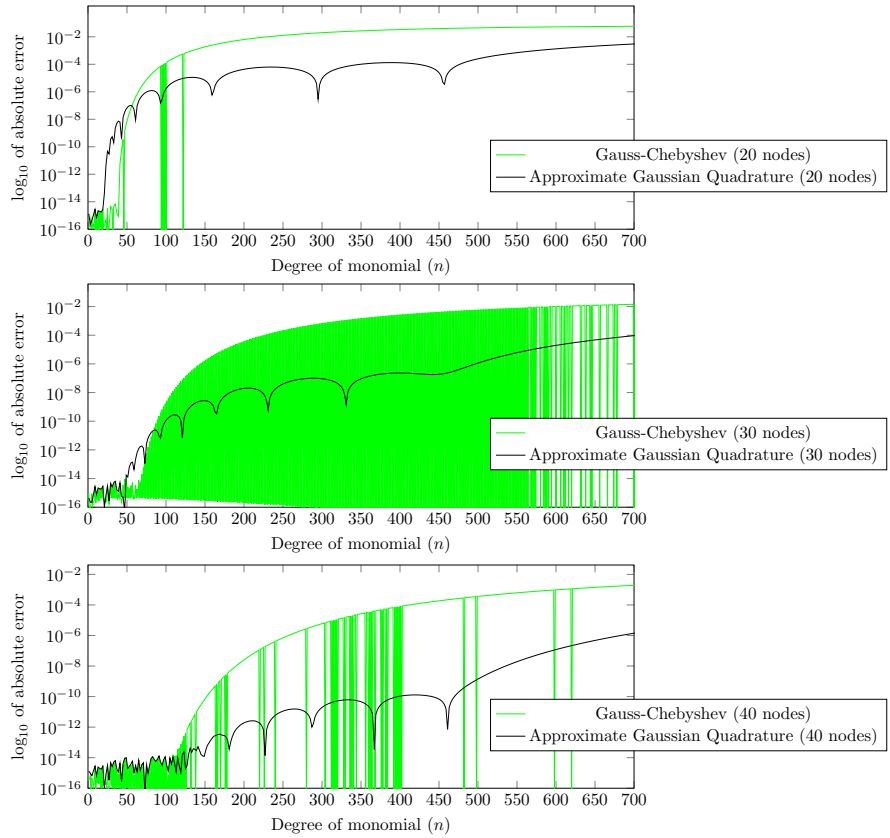

 \begin{center}
	\begin{tikzpicture}[scale=0.65]
	\pgfplotsset{every axis legend/.append style={at={(1.1,0.4)},anchor=north}}	
 	\begin{axis}[
		ylabel={$\log_{10}$ of absolute error},
	 	xlabel={Degree of monomial ($n$)},
 		ymode=log,
 		width=0.75\textwidth,
 		height=175pt,
 		xmin=0,
 		xmax=700,
 		ymin=1e-16,
 		ymax=0,
 		xtick={0, 50,100,150,200,250,300,350, 400,450,500,550,600, 650,700},
 		ytick={1e-16,1e-14, 1e-12, 1e-10,1e-08, 1e-06, 1e-4, 1e-2}]

	\input{GaussChebyshev20}
 
 	\legend{
			Gauss-Chebyshev (20 nodes) \\%
		Approximate Gaussian Quadrature (20 nodes) \\%
 		}
 	
 	\end{axis}
	\end{tikzpicture}%

	\begin{tikzpicture}[scale=0.65]
	\pgfplotsset{every axis legend/.append style={at={(1.1,0.4)},anchor=north}}	
 	\begin{axis}[
			ylabel={$\log_{10}$ of absolute error},
	 	xlabel={Degree of monomial ($n$)},
 		ymode=log,
 		width=0.75\textwidth,
 		height=175pt,
 		xmin=0,
 		xmax=700,
 		ymin=1e-16,
 		ymax=0,
 		xtick={0, 50,100,150,200,250,300,350, 400,450,500,550,600,650,700},
 		ytick={1e-16,1e-14, 1e-12, 1e-10,1e-08, 1e-06, 1e-4, 1e-2}]

	\input{GaussChebyshev30}
 
 	\legend{
		Gauss-Chebyshev (30 nodes) \\%
		Approximate Gaussian Quadrature (30 nodes) \\%
 		}
 	
 	\end{axis}
\end{tikzpicture}%

\begin{tikzpicture}[scale=0.65]
	\pgfplotsset{every axis legend/.append style={at={(1.1,0.4)},anchor=north}}	
 	\begin{axis}[
			ylabel={$\log_{10}$ of absolute error},
	 	xlabel={Degree of monomial ($n$)},
 		ymode=log,
 		width=0.75\textwidth,
 		height=175pt,
 		xmin=0,
 		xmax=700,
 		ymin=1e-16,
 		ymax=0,
 		xtick={0,50,100,150,200,250,300,350, 400,450,500,550,600,650,700},
 		ytick={1e-16,1e-14, 1e-12, 1e-10,1e-08, 1e-06, 1e-4, 1e-2}]

	\input{GaussChebyshev40}
 
 	\legend{
		Gauss-Chebyshev (40 nodes) \\%
		Approximate Gaussian Quadrature (40 nodes) \\%
 		}
 	
 	\end{axis}
	\end{tikzpicture}%

	\end{center}
 	\caption{Comparison between the absolute error incurred in the evaluation of the integral $\int_{-1}^1 x^n \, \frac{1}{\sqrt{1-x^2}} \d x$ through an Approximate Gaussian Quadrature of order $N=350$ (black)  and a Gauss-Chebyshev quadrature (green) for different number of nodes. Top: 20 nodes, Middle: 30 nodes, Bottom: 40 nodes}
	\label{GCcompare}
\end{figure}

Numerical results are shown in Figure \ref{GLcompare} and \ref{GCcompare}. They were obtained using $N=350$. The results need to be interpreted carefully. The choice of $N$ represents in effect the polynomial order that would be required to approximate a given function $f(x)$ to some accuracy $\epsilon$. A numerical quadrature will then be able to approximate the integral of $f(x)$ if it can integrate all monomials of degree less than $N$ with accuracy $\epsilon$. In Fig.~\ref{GLcompare} for example, we see that the Gauss-Legendre quadrature is exact to machine precision up to $n=39$. However the error increases rapidly to reach $10^{-3}$ near $n=350$. In contrast, although AGQ is not exact for $n \le 39$, the error up to $n \le 350$ remains lower than $10^{-4}$ with only 20 nodes. As we increase the number of nodes (middle and bottom plots) the gain below $n= 350$ is even more significant.

The behavior of AGQ in the top plot around $n \approx 40$ where Gauss-Legendre seems to outperform AGQ is not significant. Indeed if a polynomial of order $n \approx 40$ is sufficient to approximate $f$, we would reduce $N$. This would result in an AGQ quadrature much more accurate in the range $n \in [0,40]$.
 
\begin{figure}[htbp]
 \begin{center}	
	\begin{tikzpicture}[scale=0.65]
	\pgfplotsset{every axis legend/.append style={at={(1.1,0.4)},anchor=north}}	
 	\begin{axis}[
	xlabel={Degree of monomial ($n$)},
 		ymode=log,
 		width=0.75\textwidth,
 		height=175pt,
 		xmin=0,
 		xmax=350,
 		ymin=1e-17,
 		ymax=0,
 		xtick={0, 50,100,150,200,250,300,350},
 		ytick={1e-16,1e-14, 1e-12, 1e-10,1e-08, 1e-06, 1e-4, 1e-2}]

	\input{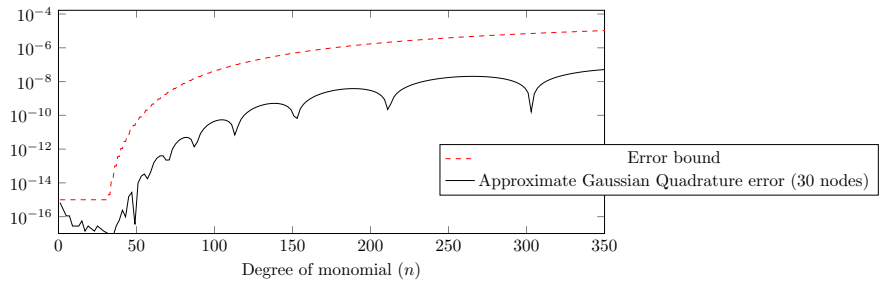}
 
 	\legend{
	 	Error bound\\%
		     Approximate Gaussian Quadrature error (30 nodes) \\%
 		}
 	
 	\end{axis}
	\end{tikzpicture}%
	\end{center}
	\caption{Comparison between the absolute error incurred in the evaluation of the integral $\int_{-1}^1 x^n \, \d x$ through a $30$-node Approximate Gaussian Quadrature of order $N=350$ and the error bound introduced in Theorem \ref{AGQ:AGQ_thm} (red)}
	\label{bound:Legendre}
\end{figure}

\begin{figure}[htbp]
	\begin{center}
	
	 \begin{tikzpicture}[scale=0.65]
	\pgfplotsset{every axis legend/.append style={at={(1.1,0.4)},anchor=north}}	
 	\begin{axis}[
	xlabel={Degree of monomial ($n$)},
 		ymode=log,
 		width=0.75\textwidth,
 		height=175pt,
 		xmin=0,
 		xmax=350,
 		ymin=1e-17,
 		ymax=0,
 		xtick={0, 50,100,150,200,250,300,350},
 		ytick={1e-16,1e-14, 1e-12, 1e-10,1e-08, 1e-06, 1e-4, 1e-2}]

	\input{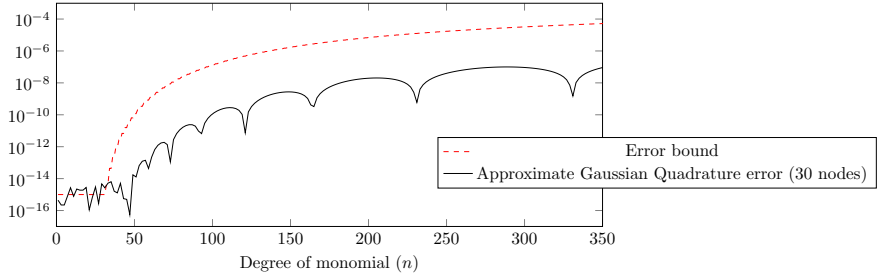}
 
 	\legend{
			Error bound\\%
		     Approximate Gaussian Quadrature error (30 nodes) \\%
 		}
 	
 	\end{axis}
	\end{tikzpicture}%

	\end{center}
		\caption{Comparison between the absolute error incurred in the evaluation of the integral $\int_{-1}^1 x^n \, \frac{1}{\sqrt{1-x^2}} \d x$ through a $30$-node Approximate Gaussian Quadrature of order $N=350$ and the error bound introduced in Theorem \ref{AGQ:AGQ_thm} (red).}
	\label{bound:Chebyshev}
\end{figure}

On Figure \ref{bound:Legendre} and \ref{bound:Chebyshev}, we also compare the theoretical bound obtained in Theorem \ref{AGQ:AGQ_thm} with the actual absolute error obtained through a 30-node AGQ for both the Lebesgue and Chebyshev measures respectively. In both cases, it is seen that the bound provides a reasonable estimate for the behavior of the error.

Finally, an interesting thing to be noted is that in both cases the nodes associated with the approximate Gaussian quadratures were \emph{real} and the weights were \emph{real and positive}; it is a known fact that this should be the case for classical Gaussian quadratures. However, this is by no means obvious for the case of approximate Gaussian quadratures, and we currently have no theory demonstrating that it is always the case for real positive measures.

\subsubsection{General integrands}

An important difference between AGQ and Gaussian quadratures is that AGQ takes $N$ as a parameter. $N$ represents in effect the order of a polynomial that can approximate $f(x)$ to the desired accuracy. This is function-dependent and therefore may need to be adjusted in AGQ depending on the integrand, if one wishes to have a near optimal quadrature.

Generally speaking, AGQ should be able to outperform a classical Gaussian quadrature in all cases since a Gaussian quadrature is a special case of AGQ, by basically choosing $d=N-1$ where $d$ is the degree of the polynomial. Indeed, this is what we observed in our numerical tests. Whenever the classical Gaussian quadrature or CGQ performs well, no gain is obtained with AGQ. We note that, in this case, the usual numerical techniques to evaluate Gaussian quadrature nodes should be more effective than the numerical procedure we are advocating for AGQ (due to ill-conditioning for too stringent a tolerance as mentioned in the introduction).


Conversely, when the convergence of CGQ is slow, AGQ provides a significant improvement. This corresponds to situation where expanding $f$ using polynomials requires terms of high degree and then the approximation of AGQ for high order monomials makes a difference. This is illustrated in the examples below.

We used the following integrands to investigate the accuracy of AGQ:
\begin{align*}
\log \left ( 1- \frac{x}{1.05} \right ) &= - \sum_{n=0}^{\infty} \frac{1}{n (1.05)^n} \,x^n ,\;\;\; |x| \leq 1 \\
\frac{1}{ 1- \frac{x}{1.05} } &= - \sum_{n=0}^{\infty} \frac{1}{(1.05)^n} \,x^n ,\;\;\; |x| \leq 1 \\
e^{-10\,x} &= - \sum_{n=0}^{\infty} \frac{(-10)^n}{n!} \,x^n ,\;\;\; 0\leq x \leq 1
\end{align*}
The first two integrand have slowly-decaying coefficients and can be approximated in the interval $[-1,1]$ through a sum containing $\O( \log_{1.05}(1/\epsilon) ) $ terms for an accuracy of $\epsilon$. At $\epsilon$-machine ($\epsilon = 10^{-15}$) this implies approximately $700$ terms. The third integrand has very fast decay, and in this case only $50$ terms are sufficient.

For each case, we varied the number of nodes in the quadrature. Then for AGQ, we selected the integer $N$ that gave us the most accurate result. In practice, an algorithm would be required to estimate $N$ numerically but we will not address this question here. Results are show in Table \ref{general:log}--\ref{general:exp}.

\begin{table}[htbp]
\begin{center}
\begin{tabular}{cccc} \toprule
Number of nodes & Optimal value for $N$ & AGQ & Gauss-Legendre\\ \midrule
10 & 75 &  $2.16 \cdot 10^{-8}$ & $1.39 \cdot 10^{-4}$ \\
15 & 100 & $1.08 \cdot 10^{-8}$ & $3.94 \cdot 10^{-6}$ \\
20 & 150 & $2.05 \cdot 10^{-11}$ & $1.26 \cdot 10^{-7}$ \\
25 & 200 & $3.99 \cdot 10^{-14}$ & $4.31 \cdot 10^{-9}$ \\
30 & 250 & $1.61 \cdot 10^{-15}$ & $1.54 \cdot 10^{-10}$ \\ \bottomrule
\end{tabular}
\end{center}
\caption{Absolute error incurred by an AGQ and a Gauss-Legendre quadrature for the integration of $f(x) = \log \left ( 1- \frac{x}{1.05}   \right ) $ over the interval $[-1,1]$ for various number of nodes.}
\label{general:log}
\end{table}%

\begin{table}[htbp]
\begin{center}
\begin{tabular}{cccc} \toprule
Number of nodes & Optimal value for $N$ & AGQ & Gauss-Legendre\\ \midrule
10 & 75 & $5.81 \cdot 10^{-5}$ & $8.15 \cdot 10^{-3}$ \\
15 & 100 & $2.20 \cdot 10^{-6}$ & $3.60 \cdot 10^{-4}$ \\
20 & 150 & $4.26 \cdot 10^{-9}$ & $1.56 \cdot 10^{-5}$ \\
25 & 200 & $1.58 \cdot 10^{-11}$ & $6.76 \cdot 10^{-7}$ \\
30 & 250 & $4.01 \cdot 10^{-13}$ & $2.92 \cdot 10^{-8}$ \\
35 & 300 & $1.77 \cdot 10^{-15}$ & $1.25 \cdot 10^{-9}$ \\ \bottomrule
\end{tabular}
\end{center}
\caption{Absolute error incurred by an AGQ and a Gauss-Legendre quadrature for the integration of $f(x) = \frac{1}{ 1- \frac{x}{1.05} }$ over the interval $[-1,1]$ for various number of nodes.}
\label{general:geo}
\end{table}%

\begin{table}[htbp]
\begin{center}
\begin{tabular}{cccc} \toprule
Number of nodes & Optimal value for $N$ & AGQ & Gauss-Legendre\\ \midrule
5 & 15 &  $1.09 \cdot 10^{-6}$ & $8.82 \cdot 10^{-5}$ \\
7 &  7 & $1.29 \cdot 10^{-7}$ & $1.29 \cdot 10^{-7}$ \\
10 & 10 & $1.02 \cdot 10^{-12}$ & $1.02 \cdot 10^{-12}$ \\
12 & 12 & $4.44 \cdot 10^{-16}$ & $4.44 \cdot 10^{-16}$ \\ \bottomrule
\end{tabular}
\end{center}
\caption{Absolute error incurred for $e^{-10\,x} $ over the interval $[0,1]$. In that case, Gauss-Legendre converges very fast and AGQ simply provides a quadrature with the same accuracy. The two methods become essentially identical.}
\label{general:exp}
\end{table}

We observe the superior accuracy of AGQ. The first two cases are challenging for CGQ and AGQ does significantly better. For the last case, CGQ converges extremely fast and then AGQ simply finds that the optimal choice is CGQ and provides an estimate with the same accuracy.

In summary, $N$ shoud be adjusted depending on the type of integrand. If the integrand is such that expansions in a polynomial basis possess slowy-decaying coefficients, AGQ will provide significantly greater accuracy. If on the contrary, a polynomial expansion converges very rapidly, both AGQ and CGQ will provide essentially identical (and fast) convergence.

We also stress that AGQ can be constructed for a wide range of measures whereas CGQ is restricted to positive measures (weight function) only.

\subsection{Singular functions}
\label{NS:sing}

\begin{figure}[htbp]
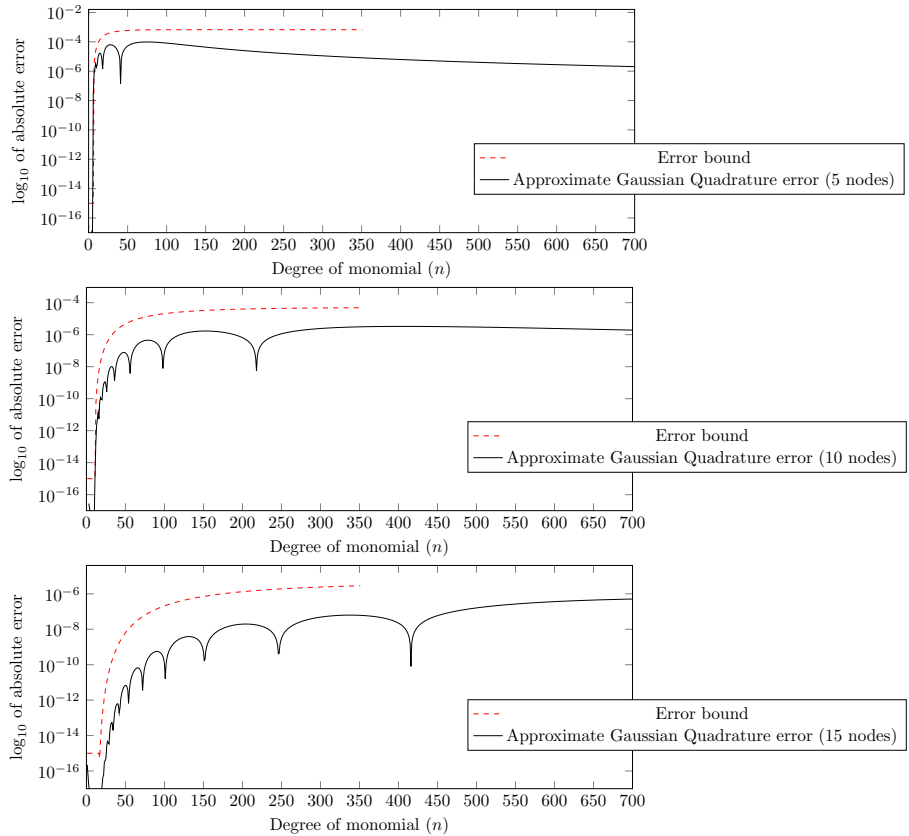

\begin{center}
	
\begin{tikzpicture}[scale=0.65]
	\pgfplotsset{every axis legend/.append style={at={(1.1,0.4)},anchor=north}}	
 	\begin{axis}[
				ylabel={$\log_{10}$ of absolute error},
	 	xlabel={Degree of monomial ($n$)},
 		ymode=log,
 		width=0.75\textwidth,
 		height=175pt,
 		xmin=0,
 		xmax=700,
 		ymin=1e-17,
 		ymax=0,
 		xtick={0, 50,100,150,200,250,300,350,400,450,500,550,600,650,700},
 		ytick={1e-16,1e-14, 1e-12, 1e-10,1e-08, 1e-06, 1e-4, 1e-2}]

	\input{LogMeasure_5}
 
 	\legend{
			Error bound\\%
		     Approximate Gaussian Quadrature error (5 nodes) \\%
 		}
 	
 	\end{axis}
\end{tikzpicture}%
	
\begin{tikzpicture}[scale=0.65]
	\pgfplotsset{every axis legend/.append style={at={(1.1,0.4)},anchor=north}}	
 	\begin{axis}[
				ylabel={$\log_{10}$ of absolute error},
	 	xlabel={Degree of monomial ($n$)},
 		ymode=log,
 		width=0.75\textwidth,
 		height=175pt,
 		xmin=0,
 		xmax=700,
 		ymin=1e-17,
 		ymax=0,
 		xtick={0, 50,100,150,200,250,300,350,400,450,500,550,600,650,700},
 		ytick={1e-16,1e-14, 1e-12, 1e-10,1e-08, 1e-06, 1e-4, 1e-2}]

	\input{LogMeasure_10}
 
 	\legend{
			Error bound\\%
		    Approximate Gaussian Quadrature error (10 nodes) \\%
 		}
 	
 	\end{axis}
\end{tikzpicture}%
	
\begin{tikzpicture}[scale=0.65]
	\pgfplotsset{every axis legend/.append style={at={(1.1,0.4)},anchor=north}}	
 	\begin{axis}[
		ylabel={$\log_{10}$ of absolute error},
	 	xlabel={Degree of monomial ($n$)},
 		ymode=log,
 		width=0.75\textwidth,
 		height=175pt,
 		xmin=0,
 		xmax=700,
 		ymin=1e-17,
 		ymax=0,
 		xtick={0, 50,100,150,200,250,300,350,400,450,500,550,600,650,700},
 		ytick={1e-16,1e-14, 1e-12, 1e-10,1e-08, 1e-06, 1e-4, 1e-2}]

	\input{LogMeasure_15}
 
 	\legend{
			Error bound\\%
		     	Approximate Gaussian Quadrature error (15 nodes) \\%
 		}
 	
 	\end{axis}
	\end{tikzpicture}%

	\end{center}
	 	\caption{Comparison between the absolute error incurred in the evaluation of the integral $\int_{-1}^1 x^n \, \log(x)  \d x$ through an Approximate Gaussian Quadrature of order $N=350$ (black)  and the error bound introduced in Theorem \ref{AGQ:AGQ_thm} (red) for 5, 10 and 15 nodes. \label{NS:log}}
 \end{figure}

We show how AGQ can be used to integrate functions with integrable singularities. For this purpose, we consider integrand of the form $x^n \log(x)$ for $x \in (0,1]$ and $0 \leq n \leq 700$. In this case, the integral of interest takes the form, 
\begin{equation*}
\int_{0}^1 x^n \, \log(x) \, \d x
\end{equation*}
This quantity can either be seen as the integration of $x^n \, \log(x)$ with respect to Lebesgue measure or as the integration of the monomial $x^n$ with respect to the measure $\d \alpha(x) =  \log(x) \, \d x$. Considering the latter, we build an AGQ of order $N=350$ with different number of nodes and display the absolute error as a function of the degree $n$ and the number of quadrature points. This is shown in Figure \ref{NS:log}. Note that the bound is not plotted beyond $N=350$ for it is no more valid past this point.

We note that in  this case we cannot perform a comparison with a classical Gaussian quadrature for no such quadrature exists as is the case with most measures but a few.

\subsection{Quadrature for polynomials on the complex circle}
\label{NS:trig}
In this section, we are interested in integrands that take the form of trigonometric monomials, i.e., functions of the form,
\begin{equation*}
f(x) =  e^{\i n x}
\end{equation*}
where $0 \leq n$. As their name conveys, such functions are just homogeneous polynomials $z^n$ in the complex plane which have been restricted to the boundary of the unit circle, i.e., $z = e^{ix}$. Thanks to this close relationship with polynomials on the real axis, one can also develop approximate Gaussian quadratures for such functions as well. In fact it suffices to replace the moments $\mu_n$ by the trigonometric moments,
\begin{equation*}
\tau_n = \int (e^{\i x})^n \, \d \alpha(x)
= \int z^n \, \d \alpha(z)
\end{equation*}
in all that has been presented above and similar results follow.

\begin{figure}[htbp]
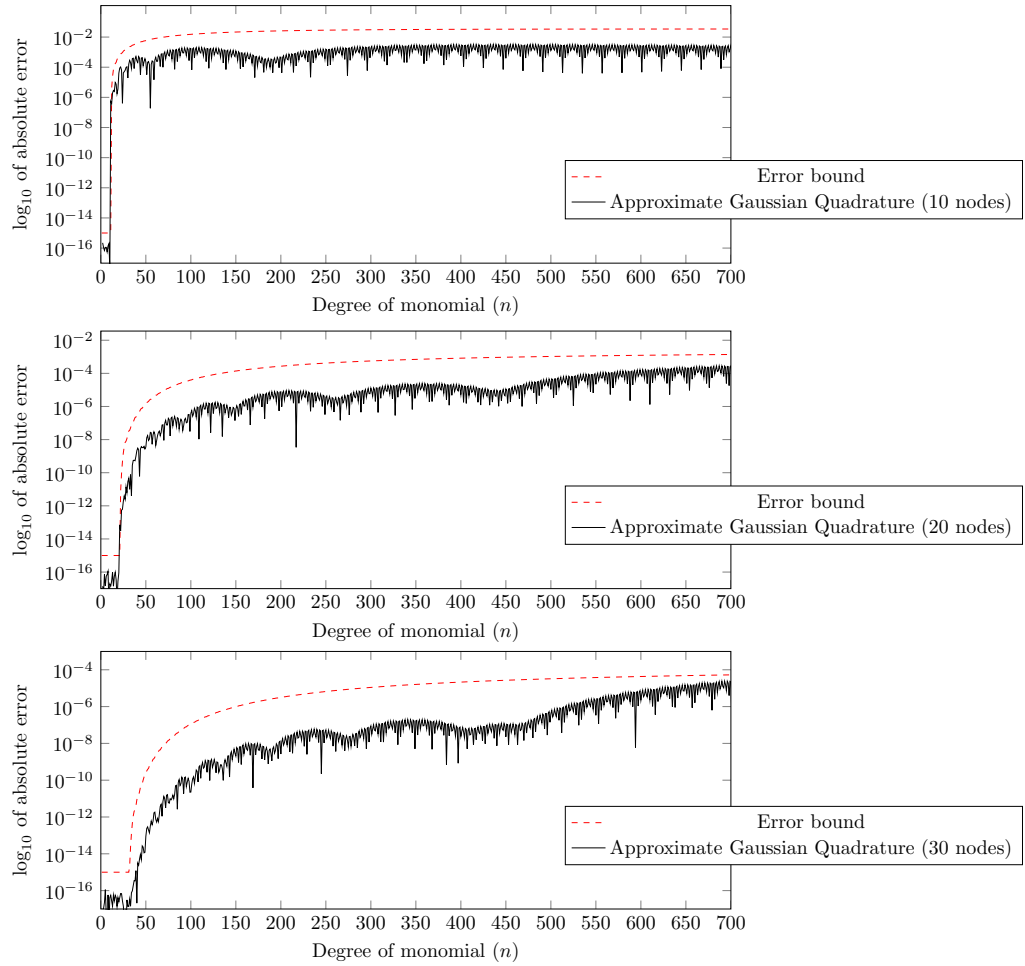

 \begin{center}
	 \begin{tikzpicture}[scale=0.75]
	 \pgfplotsset{every axis legend/.append style={at={(1.1,0.4)},anchor=north}}	
 	\begin{axis}[
		ylabel={$\log_{10}$ of absolute error},
	 	xlabel={Degree of monomial ($n$)},
 		ymode=log,
 		width=0.75\textwidth,
 		height=175pt,
 		xmin=0,
 		xmax=700,
 		ymin=1e-17,
 		ymax=0,
 		xtick={0,50,100,150,200,250,300,350,400,450,500,550,600,650,700,750,800},
 		ytick={1e-16,1e-14, 1e-12, 1e-10,1e-08, 1e-06, 1e-4, 1e-2}]
		
	\input{OscillatoryIntegrand_10}
 
 	\legend{
			Error bound\\%
		Approximate Gaussian Quadrature (10 nodes) \\%
 		}
 	
 	\end{axis}
\end{tikzpicture}%
	
\begin{tikzpicture}[scale=0.75]
	\pgfplotsset{every axis legend/.append style={at={(1.1,0.4)},anchor=north}}	
 	\begin{axis}[
		ylabel={$\log_{10}$ of absolute error},
	 	xlabel={Degree of monomial ($n$)},
 		ymode=log,
 		width=0.75\textwidth,
 		height=175pt,
 		xmin=0,
 		xmax=700,
 		ymin=1e-17,
 		ymax=0,
 		xtick={0,50,100,150,200,250,300,350,400,450,500,550,600,650,700,750,800},
 		ytick={1e-16,1e-14, 1e-12, 1e-10,1e-08, 1e-06, 1e-4, 1e-2}]

	\input{OscillatoryIntegrand_20}
 
 	\legend{
		Error bound\\%
		Approximate Gaussian Quadrature (20 nodes) \\%
 		}
 	
 	\end{axis}
\end{tikzpicture}%
	
\begin{tikzpicture}[scale=0.75]
	\pgfplotsset{every axis legend/.append style={at={(1.1,0.4)},anchor=north}}	
 	\begin{axis}[
		ylabel={$\log_{10}$ of absolute error},
	 	xlabel={Degree of monomial ($n$)},
 		ymode=log,
 		width=0.75\textwidth,
 		height=175pt,
 		xmin=0,
 		xmax=700,
 		ymin=1e-17,
 		ymax=0,
 		xtick={0,50,100,150,200,250,300,350,400,450,500,550,600,650,700,750,800},
 		ytick={1e-16,1e-14, 1e-12, 1e-10,1e-08, 1e-06, 1e-4, 1e-2}]

	\input{OscillatoryIntegrand_30}
 
 	\legend{
			Error bound\\%
		Approximate Gaussian Quadrature (30 nodes) \\%
 		}
 	
 	\end{axis}
\end{tikzpicture}%
\end{center}
\caption{Comparison between the absolute error incurred in the evaluation of the integral $\int_{-1}^1 e^{\i n x} \, \d x$ through an Approximate Gaussian Quadrature of order $N=350$ (black) and the error bound introduced in Theorem \ref{AGQ:AGQ_thm} (red) for various number of nodes. Top: 10 nodes, Middle: 20 nodes, Bottom: 30 nodes. \label{NS:trigplot}}
\end{figure}

As an example, we built an AGQ of order $N=350$ for trigonometric polynomials with respect to the Lebesgue measure over the interval $[-1,1]$. The absolute error between our approximation and the exact value of the integral,
\begin{equation*}
\int_{-1}^1  e^{\i n x} \, \d x = \frac{e^{\i n} - e^{-\i n}}{\i n}
\end{equation*}
are presented in Figure \ref{NS:trigplot}. There, it is seen that as little as $30$ quadrature points are necessary to integrate a complex exponential with frequency $n=500$ with $\approx 10^{-6}$ accuracy. We also plotted the theoretical bound of Theorem \ref{AGQ:AGQ_thm}. Again, it appears to be a good estimate.

\begin{figure}[htbp]
\begin{center}
\includegraphics[width=10cm]{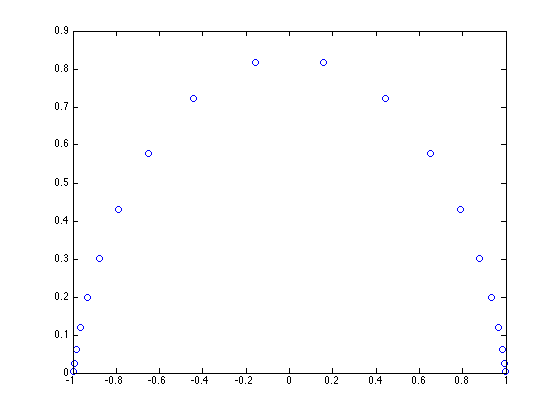}
\caption{Location in the complex plane of the nodes of a 20-node AGQ for trigonometric polynomials. The nodes appear to lie on a smooth curve with positive imaginary part.}
\label{NS:trignodes}
\end{center}
\end{figure}

It is interesting to look at the location of the nodes for such quadratures. An example is displayed in Figure \ref{NS:trignodes}. The nodes are shown in the complex plane and appear to lie along a curve which rapidly moves  upward from $-1$, slowly moves across, and rapidly moves back to $1$ on the real axis. This does not appear to be a coincidence given the fact that functions of the form $e^{\i n x}$ decay exponentially and do not oscillate along the positive imaginary axis. Thus, the underlying curve could be some sort of path of \emph{least oscillation} in an average sense over $0 \leq n \leq N$. At this point, this is a mere qualitative observation, but might be worth investigating in the future.

\subsection{Approximation of functions through short exponential sums}
\label{NS:Beylkin}
In this section, we are interested in the approximation of functions by a short sum of exponentials. That is, given a function $f(x)$ defined over an interval $[a,b]$, we seek some approximation in the form,
\begin{equation*}
f(x) \approx \sum_{m=0}^d \alpha_m e^{ \beta_m x}
\end{equation*}
for $x \in [a,b]$, and where $d$ should be as small as possible. Such expansions can be viewed as more efficient representations of functions compared to Fourier transforms as they typically require fewer terms. They can form the starting point for various fast algorithms such as the fast multipole method, hierarchical matrices ($\mathcal H$-matrices), etc. Such techniques are particularly desirable when it comes to the solution of integral equations with translation-invariant kernels (see e.g., \cite{Letourneau:2012,Beylkin:2005}). Very powerful techniques based on dynamical systems and recursion ideas were recently introduced by Beylkin \& Monz\'on \cite{Beylkin:2005, Beylkin:2010}  in order to approach this problem. As was mentioned earlier, the latter inspired the current work.

We will show how AGQ can be used to derive similar approximations through the discretization of the Fourier transform. The final formulation shares some characteristics with the problem of Beylkin \& Monz\'on that can be stated as follows: given the accuracy $\epsilon > 0$, for a smooth function $f(x)$ find the minimal number of complex weights $w_n$ and nodes $e^{t_m}$ such that,
 \begin{equation*}
\left | f(x) - \sum_m w_m e^{t_m x}   \right | < \epsilon 
\end{equation*}
for $x \in I$, $I$ being some interval in $\mathbb{R}$.

Their scheme is based on an important result regarding Hankel matrices. Consider a Hankel matrix $H$ associated with a sequence $h_k$ where $h_k = f(x_k)$ are uniform samples of $f$. Assume that the null space of $H$ is non-trivial and consider the polynomial whose coefficients are given by a vector in the null space of $H$. The zeros of this polynomial, $\lambda_i$, satisfy the following property (see e.g., \cite{Boley:1998}),
\begin{equation*}
h_k = \sum_{i=1}^r \lambda_i^k \,d_i 
\end{equation*}
for some $\{ d_i \}$, where $r$ is at most the number of columns of $H$. With our choice for $h_k$, one obtains,
\begin{equation*}
f(x_k) = \sum_{i=1}^r d_i \, e^{ \log( \lambda_i ) k } 
\end{equation*}
which naturally extends to an interpolation formula for $f(\cdot)$. 

In \cite{Beylkin:2005, Beylkin:2010}, the authors search for an approximate formula since in general the matrix $H$ is full rank and therefore no efficient representation, that would yield exactly $f(x_k)$, is possible. To achieve this, Beylkin et al.~\cite{Beylkin:2005, Beylkin:2010} show how $\lambda_i$ can be obtained as the roots of a polynomial whose coefficients are given as the entries of a con-eigenvector $u$, i.e., a vector such that,
\begin{equation*}
H u = \sigma \overline{u}
\end{equation*}
$\sigma$ being real and nonnegative. The error is then on the order of $\sigma$. They also show that the weights satisfy a well-conditioned Vandermonde system.
 
As will be seen, both our method and theirs involve a Hankel matrix with entries given by the uniform samples of the function to be approximated over the interval considered. However, the current approach avoids the solution of a con-eigenvalue problem altogether and allows for the direct computation of the weights rather than their computation through the solution of a Vandermonde system. Furthermore, since the quasi-orthogonal polynomial obtained through our scheme has small degree, the number of zeros that must be computed is also much smaller. This results in significant computational savings compared to the former method. 
  
The resulting error estimates for both methods are different. Indeed, in the case of \cite{Beylkin:2005} one expects the error to be bounded \emph{uniformly} by an expression on the order of the modulus of the small con-eigenvalue $\sigma$ (Theorem 2, \cite{Beylkin:2005}), and such value can be determined \emph{a priori}. In our case however, the error in \emph{not} uniform (as can be seen from the numerical examples). Furthermore, our current error estimate is \emph{a posteriori}.

To begin with, consider a function $f(x) \in \mathcal{L}^2 ( \mathbb{R})$ uniformly sampled at $x_n = a + \frac{n (b-a)}{N}$, $n = 0... (N-1)$ for some $N \in \mathbb{N}$ and $a,b \in \mathbb{R}$, and use the Fourier transform to write,
\begin{equation*}
f(x_n) = \int_{-\infty}^{\infty} e^{2 \pi \i x_n \xi} \hat{f} (\xi) \, \d \lambda( \xi)=\frac{N}{(b-a)}\, \int_{-\infty}^{\infty} e^{2 \pi \i n \zeta}  \, e^{2 \pi \i a \zeta} \hat{f} \left (\frac{N}{b-a} \zeta \right ) \, \d \lambda( \zeta)
\end{equation*}
where $\hat{f} (\xi)$ denotes the Fourier transform of $f(x)$, and  $\lambda( \cdot )$ is the Lebesgue measure. We note that $$\frac{N}{b-a} \,  e^{2 \pi \i a \zeta} \hat{f} \left (\frac{N}{b-a} \zeta \right )$$ can be seen as a Radon-Nykodym derivative of a certain measure $\alpha( \cdot)$ absolutely continuous with respect to Lebesgue measure (see \cite{Cohn}), i.e.,
 \begin{equation*}
 \frac{\d \alpha }{ \d \lambda} (\zeta) = \frac{N}{b-a} \,  e^{2 \pi \i a \zeta} \hat{f} \left (\frac{N}{b-a} \zeta \right )
 \end{equation*}
 With this measure we have,
 \begin{equation*}
f(x_n) = \int_{-\infty}^{\infty} e^{2 \pi \i n \zeta} \, \d \alpha(\zeta) , \;\; n = 0... N
\end{equation*}
which is perfectly well-suited for discretization through an approximate Gaussian quadrature as described in the previous section. To find such quadrature, we first need the trigonometric moments of the measure. These moments turn out to have a very simple form. Indeed, a quick look at their definition shows that,
 \begin{equation*}
\tau_n =  \int_{\mathbb{T}} e^{\i n \zeta} \, \d \alpha(\zeta) =  \int_{\mathbb{T}} e^{\i n \zeta} \, \left [ \frac{N}{b-a} \,  e^{2 \pi \i a \zeta} \hat{f} \left (\frac{N}{b-a} \zeta \right ) \right ] \d \lambda(\zeta) =  f \left( a + n \frac{(b-a)}{N} \right )
\end{equation*}
At this point, we note that the Hankel matrix arising from such moments is exactly the same as the one described in \cite{Beylkin:2005} as previously mentioned.

Finally, the nodes $\{ w_n \}$ can be obtained through Eq.\eqref{AGQ:AGQ_thm:weight}. In the end, we obtain
\begin{equation*}
f(x_n) \approx \sum_{m=0}^d w_m \, e^{  \i n \zeta_m } , \;\;\;  n=1,...,N
\end{equation*}
with error bounded by the expression provided in Theorem \ref{AGQ:AGQ_thm}. To obtain an approximation to $f(x)$ in all of $[a,b]$, we simply allow $\frac{n}{N}$ to vary continuously so that
\begin{equation*}
\frac{n}{N} =  \frac{x-a}{b-a}
\end{equation*}
for $x \in [a,b]$ and write,
\begin{align*}
f(x) &\approx \sum_{m=0}^d \alpha_m \, e^{  \i \beta_m x } \\
\alpha_m &= w_m \, e^{-\i \frac{a}{b-a} N \xi_m } \\
\beta_m &= \frac{1}{b-a} N \xi_m
\end{align*}

When $x$ corresponds to a sample, i.e., $x = x_n$ for some $n$, this reduces to the previous expression. However, when $x$ lies between two samples this last formula should be seen as an interpolation. We do not currently have the complete theory describing the interpolation error. However, it was observed numerically that such error is generally of the same order as that associated with the closest sample whenever the function $f(x)$ is sufficiently oversampled. Numerical examples are provided below. 

At this point, we describe an algorithm for the construction of such an approximation. The description can be found in pseudo-code in Algorithm \ref{approx_alg}.

\begin{algorithm2e}[htbp]
\caption{Pseudo-code to for the construction of a short exponential sum approximation of a function  $f( \cdot )$ in an interval. \label{approx_alg}}
Pick $N \in \mathbb{N}$ sufficiently large (beyond the Nyquist rate)\\
Compute $ \tau_n =  f \left( a + n \frac{(b-a)}{N} \right )$\\
Build the Hankel matrix $H_{i,j} = \tau_{i+j}$ for $i,j = 0 .. N$\\
Proceed as described in Algorithm \ref{poly_alg} to find $p(x)$\\
Compute $\{ x_n \}$, the nodes/zeros of $p(x)$ \\
Compute weights $w_n$ following Eq. \eqref{AGQ:AGQ_thm:weight} \\
Build approximation: $ \sum_{n} w_n \, e^{  \i \frac{(x-a)}{(b-a)} N \, \frac{\log(x_n)}{\i}  }$
\end{algorithm2e}

\begin{figure}[htbp]
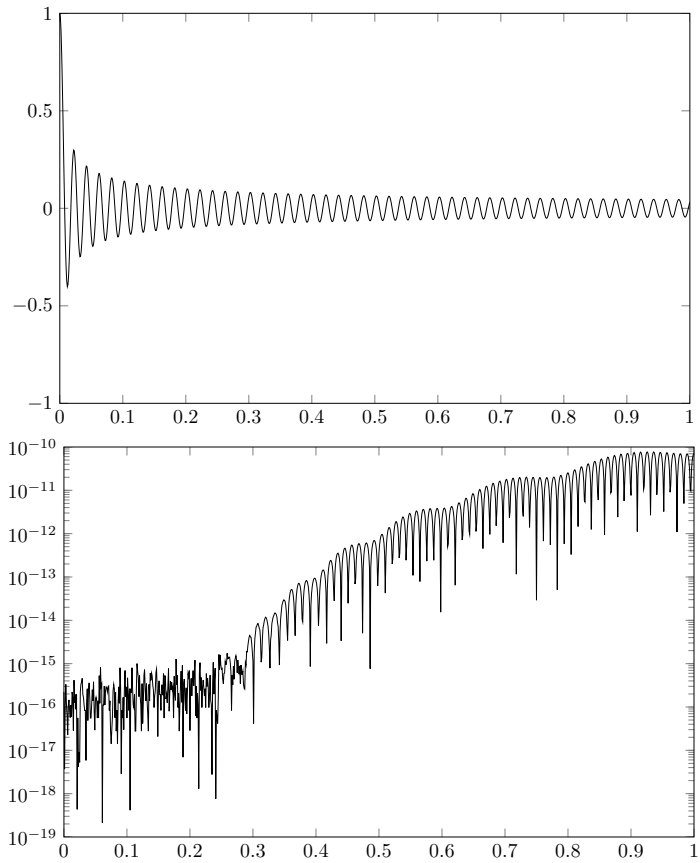

 \begin{center}
	\begin{tikzpicture}[scale=0.75]
	\pgfplotsset{every axis legend/.append style={at={(1.1,0.4)},anchor=north}}	
 	\begin{axis}[
 		width=0.75\textwidth,
 		height=0.5\textwidth,
 		xmin=0,
 		xmax=1,
 		ymin=-1.,
 		ymax=1,
 		xtick={0, 0.1, 0.2, 0.3, 0.4, 0.5, 0.6, 0.7, 0.8, 0.9, 1},
 		ytick={-1,-0.5,0,0.5,1}]

	\input{Bessel0}

 	\end{axis}
	\end{tikzpicture}%

	\begin{tikzpicture}[scale=0.75]
	\pgfplotsset{every axis legend/.append style={at={(1.1,0.4)},anchor=north}}	
 	\begin{axis}[
 		ymode=log,
 		width=0.75\textwidth,
 		height=0.5\textwidth,
 		xmin=0,
 		xmax=1,
 		ymin=1e-19,
 		ymax=1e-10,
 		xtick={0, 0.1, 0.2, 0.3, 0.4, 0.5, 0.6, 0.7, 0.8, 0.9, 1},
 		ytick={1e-19, 1e-18, 1e-17, 1e-16, 1e-15, 1e-14, 1e-13, 1e-12, 1e-11, 1e-10}]

	\input{Bessel0_error}

 	\end{axis}
	\end{tikzpicture}%

	\end{center}
	\caption{ $40$-term exponential sum approximation of Bessel function of the first kind of order $0$ ($J_{0} (100\pi x)$) in $[0,1]$ (top) and absolute error (bottom) }
	\label{bessel0}
 \end{figure}

We now provide a few examples for the representation of some oscillatory functions: the Bessel functions of the first kind $J_{\nu} (100 \pi \, x )$ over the interval $[0,1]$ and for orders $ \nu \in \{ 0, 25 \}$. Such functions are relevant in problems involving the scattering of waves in two dimensions for instance. In both cases, the order of the AGQ is $N=400$ (note that the spectrum of both functions is bounded by about $400 \approx 100 \pi$) and a $40$-terms approximation is obtained using the scheme just introduced. The results are presented in Figure \ref{bessel0} and \ref{bessel25} respectively. Agreement within $10^{-10}$ and $10^{-7}$ absolute error is observed in each cases respectively. It should also be noted that the number of terms lies much below what should be expected with a standard Fourier series given the nature of the oscillations.

 \begin{figure}[htbp]
 \begin{center}

	\begin{tikzpicture}[scale=0.75]
	\pgfplotsset{every axis legend/.append style={at={(1.1,0.4)},anchor=north}}	
 	\begin{axis}[
 		width=0.75\textwidth,
 		height=0.5\textwidth,
 		xmin=0,
 		xmax=1,
 		ymin=-0.25,
 		ymax=0.25,
 		xtick={-1, -0.8, -0.6, -0.4, -0.2, 0,  0.2,0.4,  0.6,  0.8, 1},
 		ytick={-0.3, -0.2, -0.1, 0., 0.1, 0.2, 0.3}]

	\input{Bessel25}

 	\end{axis}
	\end{tikzpicture}%
	
	\begin{tikzpicture}[scale=0.75]
	\pgfplotsset{every axis legend/.append style={at={(1.1,0.4)},anchor=north}}	
 	\begin{axis}[
 		ymode=log,
 		width=0.75\textwidth,
 		height=0.5\textwidth,
 		xmin=0,
 		xmax=1,
 		ymin=1e-17,
 		ymax=1e-7,
 		xtick={-1, -0.8, -0.6, -0.4, -0.2, 0,  0.2,0.4,  0.6,  0.8, 1},
 		ytick={1e-19, 1e-18, 1e-17, 1e-16, 1e-15, 1e-14, 1e-13, 1e-12, 1e-11, 1e-10, 1e-9, 1e-8, 1e-7, 1e-6}]

	\input{Bessel25_error}

 	\end{axis}
	\end{tikzpicture}%
	
	\end{center}
	\caption{ $40$-term exponential sum approximation of Bessel function of the first kind of order $25$ ($J_{25} (100 \pi x)$) in $[0,1]$ (top) and absolute error (bottom) }
	\label{bessel25}
\end{figure}

As a final example, we chose to represent the Dirichlet kernel,
\begin{equation*}
D_N (x) = \sum_{k=-N}^N e^{\i k x} = \frac{\sin \left (  \pi( N + 1/2) \, x  \right ) }{\sin \left (  \pi/2 \, x  \right )}
\end{equation*}
over the interval $[-1,1]$. When applied through convolution, the Dirichlet kernel acts as a low-frequency filter. In this sense, a short exponential sum approximation can be used to speed up the filtering process.

\begin{figure}[htbp]
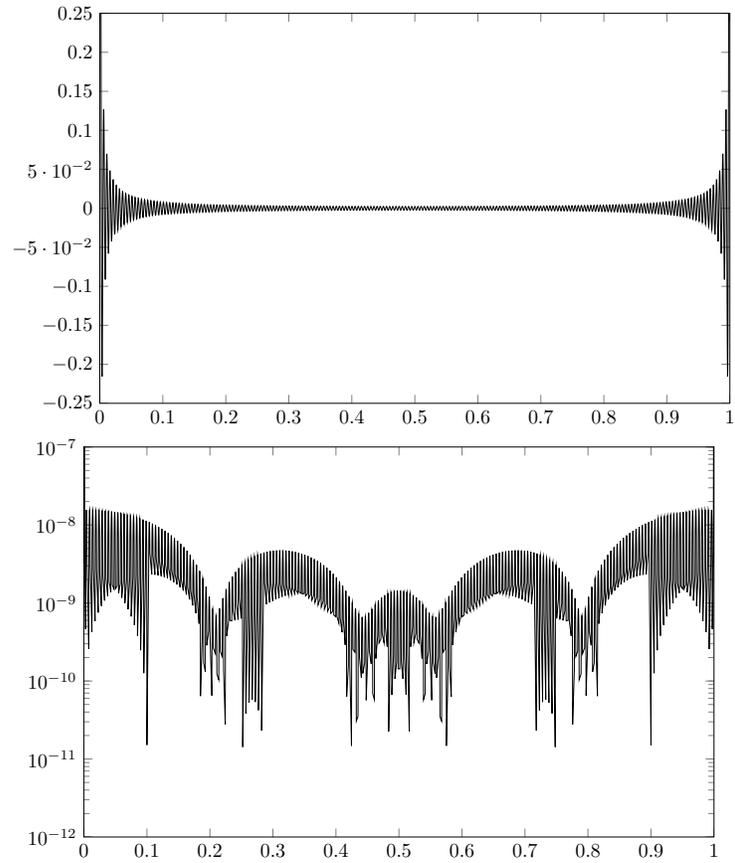

 \begin{center}
	\begin{tikzpicture}[scale=0.75]
	\pgfplotsset{every axis legend/.append style={at={(1.1,0.4)},anchor=north}}	
 	\begin{axis}[
 		width=0.75\textwidth,
 		height=0.5\textwidth,
 		xmin=0,
 		xmax=1,
 		ymin=-0.25,
 		ymax=0.25,
 		xtick={0, 0.1, 0.2, 0.3, 0.4, 0.5, 0.6, 0.7, 0.8, 0.9, 1},
 		ytick={-0.25, -0.2, -0.15, -0.1, -0.05, 0., 0.05, 0.1, 0.15, 0.2, 0.25}]

	\input{Dirichlet200}

 	\end{axis}
	\end{tikzpicture}%

	 \begin{tikzpicture}[scale=0.75]
	 \pgfplotsset{every axis legend/.append style={at={(1.1,0.4)},anchor=north}}	
 	\begin{axis}[
 		ymode=log,
 		width=0.75\textwidth,
 		height=0.5\textwidth,
 		xmin=0,
 		xmax=1,
 		ymin=1e-12,
 		ymax=1e-7,
 		xtick={0, 0.1, 0.2, 0.3, 0.4, 0.5, 0.6, 0.7, 0.8, 0.9, 1},
 		ytick={1e-19, 1e-18, 1e-17, 1e-16, 1e-15, 1e-14, 1e-13, 1e-12, 1e-11, 1e-10, 1e-9, 1e-8, 1e-7, 1e-6}]

	\input{Dirichlet200_error}

 	\end{axis}
	\end{tikzpicture}%

	\end{center}

	\caption{ $80$-term exponential sum approximation of the Dirichlet kernel of order $200$ ($D_{200} ( x)$) in $[-1,1]$ (top) and absolute error (bottom) }
	\label{dirichlet}
 \end{figure}

We picked $N = 200$. To obtain the approximation, we proceeded as described in \cite{Beylkin:2005} and went on to first approximate,
\begin{equation*}
G_{200} (x) = \sum_{k \geq 0 } \frac{\sin(200 \pi (x+k))}{200 \pi (x+k)}
\end{equation*}
through a 40-term exponential sum and then built the Dirichlet kernel through the identity,
\begin{equation*}
D_{200} (x) = G_{200} (x)  + G_{200} (1-x) 
\end{equation*}
resulting in a 80-term approximation. It is shown in Figure \ref{dirichlet}. The error is non-uniform as expected from Theorem \ref{AGQ:AGQ_thm} but still remains below $10^{-7}$ for all values in the interval.

\section{Conclusion}

We have introduced a new type of quadrature closely related to Gaussian quadratures but which use the concept of $\epsilon$-quasiorthogonality to reduce the number of quadrature nodes and weights. Such quadratures have desirable computational properties and can be applied to a family much broader than that targeted by classical Gaussian quadratures. We have provided the theory for the existence of such quadratures and have provided error estimates together with practical ways of constructing them. We have also carried out various numerical examples displaying the versatility and performance of the method. Finally, we have described how AGQ can be used to approximate functions through short exponential sums and provided further numerical examples in these cases.

\section{Acknowledgements}
The authors would like to thank Professor Ying Wu from King Abdullah University of Science and Technology (KAUST) for supporting this research through her grant as well as the National Sciences and Engineering Research Council of Canada (NSERC) for their financial support.

\newpage
\appendix

{\bf Proof of Lemma \ref{AGQ:errorlemma}.}
First, thanks to the factorization theorem for polynomials (see e.g., \cite{Hungerford})
\begin{equation}
\label{error:1}
q(x) - \tilde{q} (x) = \sum_{i=0}^{d} (q_i - \tilde{q}_i ) x^i +  \sum_{i=d+1}^{n+d} q_i x^i = \left( \sum_{i=0}^d p_i x^i \right ) \left ( \sum_{i=0}^n r_i x^i  \right ) = p(x) r(x) 
\end{equation}
and from the Cauchy product, we have,
\begin{equation}
\label{error:2}
 \left( \sum_{i=0}^d p_i x^i \right ) \left ( \sum_{i=0}^N r_i x^i  \right ) = \sum_{i=0}^{N+d} \left(  \sum_{k=0}^i r_k p_{i-k}  \right ) x^i
\end{equation}
where it is understood that coefficients corresponding to indices outside the original range of definition of $p(x)$ and $r(x)$ are $0$. By matching coefficients of like powers in Eq.\eqref{error:1} and \eqref{error:2} and putting the linear system thus obtained in matrix form, one gets
\[ \Gamma r =  \Gamma
\begin{pmatrix}
r_0  \\
r_1 \\
r_2 \\ 
\vdots \\ 
r_N 
\end{pmatrix}
=  
\begin{pmatrix}
q_0 - \tilde{q}_0  \\
\vdots \\
q_d - \tilde{q}_d \\ 
q_{d+1}\\
\vdots \\ 
q_{N+d} 
\end{pmatrix} 
= \kappa
\] 
where,
\[ \Gamma =
\begin{pmatrix}
p_0 & 0 & 0 & \cdots & 0 &0  \\
p_1 &p_0 & 0 & \cdots & 0 & 0\\
p_2 & p_1  & p_0 &  \cdots & 0 & 0\\ 
\vdots & \vdots  & \vdots &  \vdots & \vdots & \vdots \\ 
0 & 0  & 0 & \cdots & 1  & p_{d}\\
0& 0 & 0 & \cdots & 0 & 1
\end{pmatrix} \] 
$\Gamma$ is an $(N+d+1) \times (N+1)$ Toeplitz matrix characterized by the coefficients of the  known quasi-orthogonal polynomial $p(x)$. We know form the existence and uniqueness theorem for the factorization of polynomials that there exists a unique solution to the above system. We further write, (assuming $N>d$)
\[ \Gamma =
\begin{pmatrix}
\Gamma_1 \\ 
\Gamma_2
\end{pmatrix}\] 
where $\Gamma_1$ is a $d \times N$ matrix containing the first $d$ rows of $\Gamma$ and $\Gamma_2$ is a $N \times N$ matrix containing the last $N$ rows of $\Gamma$. It is to be noted that $\Gamma_2$ is an upper triangular matrix with diagonal entries all equal to $1$. Therefore, all eigenvalues of $\Gamma_2$ are equal to $1$. In particular, $\Gamma_2$ is invertible and we can write,
\[
\begin{pmatrix}
r_0  \\
r_1 \\
r_2 \\ 
\vdots \\ 
r_N 
\end{pmatrix}
= \Gamma_2^{-1}
\begin{pmatrix}
q_{d+1}\\
\vdots \\ 
q_{N+d} 
\end{pmatrix} 
\] 
where $\Gamma_2^{-1} $ is also an upper triangular Toeplitz matrix with diagonal entries all equal to $1$. Therefore,
\begin{equation*}
\sum_{n=0}^N |r_n| = \lVert r \rVert_1 = \lVert \Gamma_2^{-1} \bar{q} \rVert_1
\end{equation*}

\newpage
\bibliography{biblio}
\end{document}